\documentclass{amsart}   

 \usepackage{amsmath}
 \usepackage{amssymb}
 \usepackage{amscd} 
 \usepackage[mathscr]{eucal}
 \newtheorem{thm}{Theorem}[section]
 \newtheorem{lem}[thm]{Lemma}

 \newtheorem{cor}[thm]{Corollary}
 \newtheorem{prop}[thm]{Proposition}

 \newtheorem{defn}[thm]{Definition} 
 
 \newtheorem{prob}[thm]{Problem}

 \newcommand{\Q}{\mathbb{Q}}

  \newcommand{\tM}{\tilde{M}}

 \title[$h$--vectors]{Lower bounds for $h$--vectors of $k$--CM, 
  independence and broken circuit complexes} 
 \author{E. Swartz} 
 \address{Malott Hall \\ Cornell University \\ Ithaca, NY14850} 
 \email{ebs@math.cornell.edu}
 \thanks{Partially supported by a VIGRE postdoc under NSF grant number 9983660
to Cornell University.}

\begin{document}
 \begin{abstract}
We present a number of  lower bounds for the $h$--vectors of $k$--CM, broken circuit and independence complexes.  These lead to bounds on the coefficients of the characteristic and reliability polynomials of matroids. The main techniques are the use of 
series and parallel constructions on matroids and the short simplicial 
$h$--vector for pure complexes.
 \end{abstract}

 \maketitle
 \section{Introduction}

Based on the ideas of Whitney \cite{Wh} and Rota \cite{Ro}, the broken circuit complex of a graph was introduced by Wilf in ``What 
polynomials are chromatic?'' \cite{Wi}  Extended to matroids by Brylawski 
\cite{Bry4}, its $f$--vector corresponds to the coefficients of the 
characteristic polynomial of the matroid.  The $h$--vector encodes the 
same information in a different way.  From the point of view of matroids, 
Wilf's original question becomes, ``What are the possible $f$--vectors, or 
equivalently $h$--vectors, of broken circuit complexes of matroids?''  

Cohen--Macaulay complexes cover a wide variety of examples.  In addition to 
the broken circuit and independence complexes of matroids covered here, 
Cohen--Macaulay complexes also include all triangulations of homology balls and spheres.  
In contrast to broken circuit complexes, the possible $h$--vectors (and hence 
$f$--vectors) of Cohen--Macaulay complexes have been completely characterized (see, 
for instance, \cite[Theorem II.3.3, pg. 59]{St}).  Introduced by Baclawski, doubly Cohen--Macualay 
complexes are Cohen--Macaulay complexes which neither lose a dimension nor 
lose the Cohen--Macaulay property when any vertex is removed.  Spheres are 
doubly Cohen--Macaulay but balls are not.  More generally, a 
Cohen--Macaulay complex is 
$k$--CM if it retains its dimension and is still Cohen--Macaulay whenever 
$k-1$ or fewer vertices are removed.  In addition to the independence complexes considered below, the order complex of a geometric lattice with the top and bottom points removed is $k$--CM if every  line has at least $k$ points \cite{Ba2}. 

The $h$--vectors of independence complexes of matroids are contained in the 
intersection of $h$--vectors of broken circuit complexes and $k$--CM 
complexes.  Precisely, the cone on any independence complex is a broken 
circuit complex.  In addition, if the smallest cocircuit of the matroid has 
cardinality $k,$ then its independence complex is a $k$--CM complex.  The close connection between $h$--vectors of independence complexes of matroids and reliability problems has been studied by a number of authors.  See \cite{CC} for a recent survey.

Upper bounds on all of the above complexes have been  studied. As 
they are all Cohen--Macaulay they share a common 
absolute upper bound of 
$h_i \le \binom{n-r-1+i}{i},$ where $n$ is the number of 
vertices and $(r-1)$ is the dimension of the complex.  In addition, they 
all satisfy the relative upper bound $h_{i+1} \le h^{<i>}_i$ (see Section 
\ref{k-CM} for a definition of $h^{<i>}_i$).

Our main purpose is to analyze absolute and relative lower bounds for the 
$h$--vectors of $k$--CM, broken circuit and independence complexes.  
Section \ref{face enumeration} contains the basic facts of the 
short--simplicial $h$--vector.  The main tool for providing relative lower 
bounds is equation (\ref{del by h}). The broken circuit and independence complex of a matroid are described in section  \ref{matroids}.  Sections \ref{k-CM}, \ref{BC complexes} and \ref{ind complexes} 
contain absolute and relative lower bounds for $k$--CM, broken circuit and 
independence complexes respectively.  

Throughout the paper $\Delta$ is an $(r-1)$--dimensional simplicial 
complex with vertex set $V, |V|=n.$  The link of a vertex $v \in V$ is 
$lk_{\Delta} v,$ or just $lk \ v$ if no confusion is possible. We use 
$\Delta - v$ for the complex obtained by removing $v$ and all of the faces 
which contain $v$ from $\Delta.$  Similarly, if $A \subseteq V,$ then, $\Delta -A$ 
is the complex obtained by removing all of the vertices in $A$ and any 
faces which contain one or more of those vertices.

\section{Face enumeration} \label{face enumeration}
  The combinatorics of a simplicial complex $\Delta$ can be encoded in 
  several ways.  The most direct is to let $f_i(\Delta)$ be the number of 
  faces of cardinality $i.$  For an $(r-1)$--dimensional complex the 
  $h$--vector of $\Delta$ is the sequence $(h_0(\Delta), \dots, 
  h_r(\Delta)),$ where
\begin{equation} \label{f by h}
  h_i(\Delta) = \sum^i_{j=0} (-1)^{i-j} \binom{r-j}{r-i} f_j(\Delta).
\end{equation}  
Equivalently,
\begin{equation} \label{h by f}
  f_j(\Delta) = \sum^j_{i=0} \binom{r-i}{r-j} h_i(\Delta).
\end{equation}

\noindent By convention, $h_i(\Delta) = f_i(\Delta) =0 $ if $i <0$ or $i > 
r.$ 
The {\it short simplicial} $h$--vector was introduced in \cite{HN} as a 
simplicial analogue of the short cubical $h$--vector in \cite{Ad}. It is the 
sum of the $h$--vectors of the links of the vertices.  
As far as we know, (\ref{short by h}) was first stated in \cite{Mc}.  
However, only a proof for shellable $\Delta$ was given there. So, we 
include a proof for arbitrary pure complexes for the sake of completeness.

\begin{defn}
  Let $\Delta$ be a pure  simplicial complex .  Define
\begin{equation} \label{short h}
  \tilde{h}_i(\Delta) = \sum_{v \in V} h_i ( lk \  v).
\end{equation}
\end{defn}

\begin{lem} \cite{HN}
  Let $\Delta$ be a pure  simplicial complex.  For all 
  $i, 0 \le i \le r-1,$
  \begin{equation} \label{short by f}
    \tilde{h}_i(\Delta) = \sum^i_{j=0} (-1)^{i-j} (j+1) \binom{r-j-1}{r-i-1} 
    f_{j+1}.
  \end{equation}
\end{lem}

\begin{prop}
  Let $\Delta$ be a pure  simplicial complex.  Then,
     \begin{equation} \label{short by h}
        \tilde{h}_{i-1}(\Delta) = i \ h_i(\Delta) + (r-i+1) h_{i-1}(\Delta).
      \end{equation}
If $\dim (\Delta-v) = r-1$ for every vertex $v,$ then
      \begin{equation} \label{del by h}
        \sum_{v \in V} h_i(\Delta - v) = (n-i) h_i(\Delta) - (r-i+1) 
        h_{i-1}(\Delta).
       \end{equation}
\end{prop}

\begin{proof}
  Combining (\ref{h by f}) and (\ref{short by f}),
$$ \begin{array}{lcl}
\tilde{h}_{i-1}(\Delta) &  =  & \displaystyle\sum^{i-1}_{j=0} (-1)^{i-j-1} (j+1) 
\binom{r-j-1}{r-i} \displaystyle\sum^{j+1}_{k=0} \binom{r-k}{r-j-1} h_k(\Delta) \\
\ & = & \displaystyle\sum^i_{k=0} h_k(\Delta) \left\{ \displaystyle\sum^{i-1}_{j=k-1}
(-1)^{i-j-1} 
(j+1) \binom{r-j-1}{r-i} \binom{r-k}{r-j-1} \right\} \\
\ & = & \displaystyle\sum^i_{k=0} h_k(\Delta) 
\left\{ \displaystyle\sum^{i-1}_{j=k-1} (-1)^{i-j-1} (j+1)
\binom{r-j-1}{i-j-1} \binom{r-k}{j+1-k} \right\}.
\end{array}$$
Substituting $s=j-k+1$ and $t=i-j-1,$
$$\begin{array}{lcl}
\tilde{h}_{i-1}(\Delta) &  = & \displaystyle\sum^i_{k=0} h_k(\Delta)
\left\{ \sum_{s+t=i-k} (-1)^t 
(i-t) \binom{r+t-i}{t} \binom {r+s+t-i}{s} \right\} \\
\ & = & \displaystyle\sum^i_{k=0} h_k(\Delta)
\left\{ \sum_{s+t=i-k} (-1)^t (i-t) \frac{A}{s! t!} 
\right\}, 
\end{array}$$
where $A$ is the falling factorial $(r-k) \cdot (r-k-1) \cdots (r-i+1).$

For a fixed $i$, define $c_k$ by
$$c_k= \sum_{s+t=i-k} (-1)^t (i-t) \frac{1}{s! t!}.$$
Equation (\ref{short by h}) is equivalent to showing that $c_i = i, c_{i-1} = 1$ and 
$c_k = 0$ in all other cases. This can be seen by recognizing $c_{i-k}$ as 
the $k^{th}$ term in the generating series for 
$$(i+x)e^{-x} \cdot e^x = \left( \sum^\infty_{t=0} (-1)^t \frac{(i-t)}{t!} 
x^t \right)  
\left( \sum^\infty_{s=0} \frac{1}{s!} x^s \right).$$
 
In order to prove that  (\ref{del by h}) holds, we first notice that the 
hypothesis implies that $h_i(\Delta) = h_i(\Delta-v) + h_{i-1}(lk \  v)$ for 
every vertex $v.$  Now sum this equation over all the vertices and apply equation (\ref{short by h}).
\end{proof}

The above 
proposition makes precise the idea that, taken together, $h_{i-1}(\Delta)$ 
and $h_i(\Delta)$ measure the ``average contribution of $h_{i-1}( 
lk \ v)$ to $h_i(\Delta).$''
Another consequence of (\ref{short by h}) is that if the automorphism group of a 
pure $(r-1)$--dimensional complex
$\Delta$ is transitive, or more generally if $h_{i-1}(lk \  v)$ is 
independent of $v,$ then $n$ divides $\{i\ h_i(\Delta) + (r-i+1) 
h_{i-1}(\Delta)\}.$

\section{Broken circuit and independence complexes of matroids} \label{matroids}

We follow \cite{O} for matroid terminology.  Unless otherwise  specified, $M$ is always a rank $r$ matroid with ground set $E$ (or $E(M)$ if necessary) and $|E|=n.$  There are many equivalent ways of defining matroids.  The most convenient for us is the following.  

A {\it matroid}, $M,$ is a pair $(E, \mathcal{I}), E$ a non-empty
  finite ground set and $\mathcal{I}$ a distinguished set of subsets
  of $E.$ The members of $\mathcal{I}$ are called the {\it
  independent} subsets of $M$ and are required to satisfy:

\begin{enumerate}
  \item The empty set is in $\mathcal{I}.$ 
  
  \item If $B$ is an independent set and $A \subseteq B,$ then $A$ is
    an independent set.  
   
  \item If $A$ and $B$ are independent sets such that $|A| < |B|,$
    then there exists an element $x \in B-A$ such that $A \cup x$ is
    independent.
\end{enumerate} 
Matroid theory was introduced by Whitney \cite{W}.   The prototypical
example of a matroid is a finite subset of a vector space with the canonical independent  sets.  Another example is the cycle matroid of a graph.  Here the ground set is the edge set of the graph and a collection of edges is independent if and only if it is acyclic.

An element $e$ of a matroid is a {\it loop} if it is not contained in
any independent set.   The {\it circuits} of a matroid are its minimal
dependent sets.  Every loop of $M$ is a circuit.  A maximal
independent set is called a {\it basis}, and any element which is
contained in every basis is a {\it coloop} of the matroid.  Every
basis of $M$ has the same cardinality.  The {\it rank} of $M,$ or
$r(M),$ is that common cardinality. Similarly, the rank of a subset $A$ of $E$ is the cardinality of any maximal independent subset of $A$ and is denoted $r(A).$   The {\it deletion} of $M$ at $e$ is 
denoted $M - e.$ It is the matroid whose
finite set is $E - e$ and whose independent sets are simply those
members of $\mathcal{I}$ which do not contain $e.$  The {\it
contraction} of $M$ at $e$ is denoted $M/e.$  It is a matroid whose
ground set is also $E - e.$  If $e$ is a loop or a coloop of $M$ then
$M/e = M - e.$  Otherwise, a subset $I$ of $E-e$ is independent in
$M/e$ if and only if $I \cup e$ is independent in $M.$  Deletion and
contraction for a subset $A$ of $E$ is defined by repeatedly deleting
or contracting each element of $A.$ 
  
  The {\it dual} of $M$ is $M^\star.$ It is the matroid whose ground set
is the same as $M$ and whose bases are the complements of the bases of
$M.$  For 
example, $U_{i,j}$ is the matroid defined by $E=\{1,2,\dots,j\}$ and 
$\mathcal{I}=\{A \subseteq E: |A| \le i.\}$  So, $U_{i,j}^\star = U_{j-i,j}.$
  
  Two non-loop elements $e,f \in E$ are {\it parallel} if they form a circuit.  The relation ``is parallel to'' is an equivalence relation on $E$ and the corresponding equivalence classes are the parallel classes of $M.$  If $P$ is a parallel class of $M$, then for any $e \in P$ all of the members of $P-e$ are loops in $M/e.$  A parallel class in $M^\star$ is a {\it series class} of $M.$  If $S$ is a series class of $M,$ then for any $e \in S,$ all of the members of $S - e$ are coloops in $M-e.$  

Let $M=(E,\mathcal{I})$ and $M^\prime=(E^\prime,\mathcal{I}^\prime)$
be two matroids with $E \cap E^\prime = \emptyset.$ Then $M \oplus
M^\prime$ is the  direct sum of $M$ and $M^\prime.$ It is the matroid
whose ground set is $E \cup E^\prime$ and whose independent sets are
those subsets of the form $I \cup I^\prime, I \in \mathcal{I},
I^\prime \in \mathcal{I}^\prime.$  A matroid is {\it connected} if it is not the direct sum of two smaller matroids.  Every matroid can be written uniquely (up to order) as a direct sum $M = M_1 \oplus \dots \oplus M_k$ of connected matroids.  The {\it components} of $M$ are the summands of this decomposition.
  
The {\it independence complex} of $M$ is 
$$\Delta(M) = \{A \subseteq E: A \mbox{ is independent}\}.$$
Evidently, $\Delta(M)$ is a pure 
$(r-1)$--dimensional complex, where $r$ is the rank of $M.$ In 
addition,  $\Delta(M-e) = \Delta(M) -e$ and if $e$ is not a loop of $M,$ 
then $\Delta(M/e) = lk_{\Delta(M)} \  
e.$

In order to define the broken 
circuit complex for $M$,
we first choose a linear order $\omega$ on the elements of the matroid.  
Given such 
an order, a {\it broken circuit} is a circuit with its least element 
removed. 
The {\it broken circuit complex} is the simplicial complex whose simplices 
are the subsets of $E$ which do not contain a broken circuit.  
We denote the broken circuit complex of $M$ and $\omega$ by 
$\Delta^{BC}(M),$ or $\Delta^{BC}(M,\omega).$  
Different orderings may lead to different complexes, see \cite[Example 
7.4.4]{Bj}. However, $f_i(\Delta^{BC}(M,\omega))$ does not depend on 
$\omega$ (see Theorem \ref{h,b and tp} below). 
Conversely, distinct matroids can have the same broken circuit 
complex.  For instance, let $E=\{e_1,e_2,e_3,e_4,e_5,e_6\},$ and let 
$\omega$ be the obvious order.  Let $M_1$ be the matroid on $E$ whose 
bases are all triples except $\{e_1,e_2,e_3\}$ 
and  $\{e_4,e_5,e_6\}$ and let $M_2$ be the matroid on $E$ whose 
bases are all triples except $\{e_1,e_2,e_3\}$ and 
$\{e_1,e_5,e_6\}.$  Then $M_1$ and $M_2$ are non--isomorphic matroids but 
their broken circuit complexes are identical. 

In order to easily distinguish the  $h$--vectors of $\Delta(M)$ and $\Delta^{BC}(M)$  we use the following notation.

\begin{defn} Let $M$ be a rank $r$ matroid. \label{h,b,w defs}
  \begin{itemize}
     \item
        $h_i(M) = h_i(\Delta(M)).$
     \item
        $b_i(M) = h_{r-i}(\Delta^{BC}(M)).$
     \item
        $w_i(M) = f_{r-i}(\Delta^{BC}(M)).$
     \item
        $b^\star_i(M) = b_i(M^\star) = h_{n-r-i}(\Delta^{BC}(M^\star))$
   \end{itemize}
\end{defn}

We will suppress the $M$ when there is no danger of confusion.  
The invariants $h_i, b_i, w_i, b^\star_i$ are closely related to the 
Tutte polynomial of $M.$  The {\it Tutte polynomial} is a two--variable 
polynomial invariant of $M$ defined by

$$T(M;x,y) = \displaystyle\sum_{A \subseteq E} (x-1)^{r(M)-r(A)} 
(y-1)^{|A|-r(A)}.$$

\begin{thm} \cite{Bj} \label{h,b and tp}
 Suppose $M$ has $k$ components and $j$ coloops.   Then,
\begin{itemize}
  \item[a.]
     $T(M;x,1) = h_0 x^r + h_1 x^{r-1} + \dots + h_{r-j} x^j.$
  \item[b.]
     $T(M;x,0) = b_r x^r + b_{r-1} x^{r-1} + \dots + b_k x^k.$
  \item[c.]
     $T(M;0,y) = b^\star_{n-r} x^{n-r} + \dots + b^\star_k x^k.$ 
  \item[d.]
     $(-1)^r T(M;1-x,0) = w_0 x^r - w_1 x^{r-1} + \dots  + (-1)^r w_r.$
\end{itemize}
\end{thm}

The $w_i$ are the unsigned Whitney numbers of the first kind.  The {\it 
characteristic polynomial} of $M$ is  $(-1)^r T(M;1-x,0).$
The characteristic polynomial of a matroid has a 
number of applications including graph coloring and flows, linear coding 
theory and hyperplane arrangements.   See \cite{BO} for a survey.

 Properties [a]-[d] of $b_i$ and $h_i$ listed below  follow
immediately from corresponding properties of the Tutte polynomial which can 
be found in \cite{Bry2}.  The parallel and series connection of two (pointed)  matroids is described in \cite[Section 7.1]{O}.

\begin{thm}[Tutte recursion] \ \label{tuttepoly}

\begin{itemize}
  \item[a.]
     If $M$ has $j$ coloops, then $h_i (M) = h_i (\tM),$ where $\tM$ is $M$ with the coloops deleted.  In particular, $h_i(M) > 0$ if and  only if 
     $0 \le i \le r-j.$
  \item[b.]
     If $M$ has $k$ components and no loops, then $b_i > 0$ 
     if and only if $ k \le i \le r.$  
   \item[c.]
     If $e$ is neither a loop nor a coloop of $M,$ then $h_i(M) = 
     h_i(M-e) + h_{i-1}(M/e)$ and $b_i(M) = b_i(M-e) + b_i(M/e).$
   \item[d.]
     If $M = M_1 \oplus M_2,$ then $h_i(M) = \displaystyle\sum_{j+k=i} 
     h_j(M_1) h_k(M_2)$ and \\ $b_i(M) = \displaystyle\sum_{j+k = i} b_j(M_1) 
     b_k(M_2).$
   \item[e.]
     Suppose that $P$ is a parallel class of $M.$  Let $\tM$ 
     be $M$ with all but one element, say $e,$ of $P$ deleted. Then, $h_i(M) = h_i(\tM) + (|P|-1) h_{i-1}(\tM/e).$
   \item[f.]
     Let $S$ be a series class of $M.$  Let  $\tM$ be $M$ with all but one element, say $e,$ of $S$ 
     contracted. Then $b_i(M) = b_i(\tM) + \sum^{|S|-1}_{j=1} b_{i-j}(\tM-e).$
   \item[g.]
     Let $M$ be a parallel connection of $A$ and $B,$ where the rank of 
     $A$ is $r(A)$ and the rank of $B$ is $r(B).$  The rank of $M$ is 
     $r(A) + r(B) -1.$ In addition,  $b_i(M) = \displaystyle\sum_{j+k=i+1} b_j(A) b_k(B).$ If $A$ and $B$ are connected, then $M$ is also 
     connected.  
\end{itemize}
\end{thm}
  
  \begin{proof}
    Property [g] follows from the fact that if $M$ is a parallel connection of $A $ and $B,$ then $T(M;x,0) = T(A;x,0) * T(B;x,0)/x$ \cite[pg. 179--182]{Bry2}.  Both [e] and [f] are proved by  deleting and contracting all the elements of the given parallel or series class except $e.$     
   \end{proof}
   
   One of the consequences of [a] and [f] above is that if we increase the size of a series class of cardinality $k$ in $M$ by one, then $b_1,\dots,b_k$ are unchanged, while $b_i$ for $i > k$ may increase.    

\section{Cohen--Macaulay and $k$--CM complexes} \label{k-CM}

There are several equivalent definitions of Cohen--Macaulay complexes.  
The following will suffice for our purposes.

\begin{defn}
  A pure $(r-1)$--dimensional complex $\Delta$ is {\it Cohen--Macaulay} if 
  for  every face $F \in \Delta$ and $i < \dim (lk \  F), 
  \tilde{H}_i(lk \  F;\Q)=0.$
\end{defn}   

A numerical description of all possible $h$--vectors of Cohen--Macaulay 
complexes can be given using the following operator.  Given  any positive integers $h$ and $i$
there is a 
unique way of writing
$$h=\binom{a_i}{i} + \binom{a_{i-1}}{i-1} + \dots + \binom{a_j}{j}$$
so that $a_i > a_{i-1} > \dots > a_j \ge j \ge 1.$  Define
$$h^{<i>} = \binom{a_i+1}{i+1} + \binom{a_{i-1}+1}{i} + \dots + 
\binom{a_j+1}{j+1}$$

\begin{thm} \cite{St} \label{CM h-vectors}
  A sequence of non--negative integers $(h_0, \dots, h_r)$ is the 
  $h$--vector of some Cohen--Macaulay complex if and only if $h_0=1$ and 
  $h_{i+1} \le h^{<i>}_i$ for all $1 \le i \le r-1.$
\end{thm}

The notion of $k$--CM complexes was introduced by Baclawski \cite{Ba2}.

\begin{defn}
  Let $\Delta$ be a pure $(r-1)$--dimensional simplicial complex with 
  vertex set $V$ and $k \ge 1.$  We say 
  that $\Delta$ is $k$--CM if for all $A \subseteq V$ with $|A|<k, \Delta -A$ is Cohen--Macaulay of dimension $(r-1).$
\end{defn}

Examples of $2$--CM complexes include order complexes of geometric lattices, finite buildings and triangulations of spheres.  Several examples and constructions 
involving $k$--CM complexes, especially for order complexes of posets, are 
contained in \cite{Ba2}.  Since $lk_{\Delta} v - A = lk_{\Delta - A} v,$ 
the link of any vertex of a $k$--CM complex is $k$--CM, and removing a 
vertex from a $k$--CM complex leaves a $(k-1)$--CM complex (as long as $k > 
1$).

The independence and broken circuit complexes of a matroid are Cohen--Macaulay \cite{St1}. 
So, $\Delta(M)$ is $k$--CM if and only if every hyperplane of $M$ has 
cardinality at most $n-k$. Equivalently, the smallest cocircuit of $M$ has at least $k$ elements. However, $\Delta^{BC}(M)$ is a cone on the least 
element, hence it is only $1$--CM.  If the cone point is removed, then  the 
remaining complex is also Cohen--Macaulay, but may still be only $1$--CM.  
For example, let $M$ be the cycle matroid of the theta--graph with three paths each of length 2.  Direct computation shows that the $h$--vector of 
$\Delta^{BC}(M)$ is $(1,2,3,1).$ Removing the cone point leaves a 
$2$--dimensional complex with $5$ points and
the same $h$--vector. By Corollary \ref{relative k-cm} below, $(1,2,3,1)$ is not 
the $h$--vector of any $2$--dimensional $2$--CM complex with $5$ points.

Theorem \ref{CM h-vectors} gives an upper bound for possible $h$--vectors of 
Cohen--Macaulay complexes.  It also makes it clear that there are no lower 
bounds.  For $k$--CM complexes we have the following absolute lower bound.  Recall that $U_{r,n}$ is the rank $r$ matroid with $n$ elements such that every $r$--element subset is a basis. 

\begin{prop} \label{k-CM abs}
  Let $\Delta$ be an $(r-1)$--dimensional $k$--CM complex.  
  Then,
   $$h_i(\Delta) \ge h_i(U_{r,r+k-1}).$$ 
\end{prop}

\begin{proof}
  Induction on $n$ and $k$.  When $k=1,$ the theorem is simply the 
  statement that $h_i(\Delta) \ge 0$ for $i \ge 1,$ and $h_0(\Delta) \ge 
  1.$  For fixed $k,$ the definition of $k$--CM forces $n \ge r+k-1.$  Suppose 
  $n=r+k-1.$  Since the removal of any subset of vertices of cardinality 
  $k-1$ does not lower the dimension of $\Delta,$ every subset of vertices 
  of cardinality $r$ must be a face of $\Delta.$  So, $\Delta = 
  \Delta(U_{r,r+k-1}).$    For the induction step, let $v$ 
  be any vertex of $\Delta.$  Then
$$h_i(\Delta) = h_i(\Delta -v) + h_{i-1}(lk_{\Delta} v) \ge 
h_i(U_{r,r+k-2}) + h_{i-1}(U_{r-1,r+k-2}) = h_i(U_{r,r+k-1}).$$
\end{proof}

Minimizing $h$--vectors is closely related to the problem of finding the least reliable graph.  Let $G$ be a connected graph with $r+1$ vertices and $n$ edges.  Thus $M(G),$ the cycle matroid of $G,$ has rank $r$ and cardinality $n.$  Suppose that each edge of $G$ has equiprobability $p, 0 < p < 1$ of being deleted.  Then the probability that $G$ remains connected is $R_G(p) = (1-p)^r [h_0(M(G)^\star) + h_1 (M(G)^\star) p + \dots + h_{n-r}(M(G)^\star) p^{n-r}].$  Boesch, Satyanarayana and Suffel posed the problem of finding the minimum of $R_G(p)$ among all connected simple graphs with $r+1$ vertices and $n$ edges.  They also conjectured that a particular graph, which they called $L(r+1,n),$ would attain that lower bound \cite{BSS}.  Brown, Colbourn and Devitt further conjectured that the $h$--vector of $L(r+1,n)$ would be an absolute lower bound for  the $h$--vector of $M(G)^\star$ among all connected simple graphs with $r+1$ vertices and $n$ edges \cite{BCD}.   The original conjecture of Boesch et. al. was confirmed for $n$ greater than $\binom{r-1}{2}$ in \cite{PSS}.  The corresponding problem in the category of matroids is to find among all rank $r$ cosimple matroids of cardinality $n$ one which minimizes the $h$--vector.  Since $M$ is cosimple if and only if $\Delta(M)$ is $3$--CM, the above proposition shows that $U_{0,n-r-2} \oplus U_{r,r+2}$ is the solution to this problem.

Combining the above proposition with (\ref{del by h}) immediately gives a 
relative lower bound.

\begin{cor} \label{relative k-cm}
  Let $\Delta$ be an $(r-1)$--dimensional $k$--CM complex with $n$ 
  vertices.  Then,
$$(n-i) h_i \ge (r-i+1) h_{i-1} + n \binom{i+k-3}{i}.$$
\end{cor}

\begin{proof}  For every vertex $v, \Delta - v$ is $(k-1)$--CM.  Now 
combine (\ref{del by h}), Proposition \ref{k-CM abs} and the fact that 
$h_i(U_{r,r+k-2}) = \binom{i+k-3}{i}.$
\end{proof}  

\begin{prob}
  Given $r,n,k$ and $i,$ what is the minimum of $h_i(\Delta)$ over all 
  $(r-1)$--dimensional $k$--CM complexes with $n$ vertices?  Does there 
  exist a $\Delta$ which attains these values?
\end{prob}

Conjecture II.6.2 in \cite{St} would imply that for $2$--CM complexes with 
$n$ equal to $r+2, 
h_i(\Delta) \ge h_i(\Delta(U_{1,2} \oplus U_{r-1,r})).$
In section \ref{ind complexes} we will give an answer to this problem for 
independence complexes of matroids when $n$ is sufficiently large.

\section{Broken circuit complexes} \label{BC complexes}

In this section we assume that $M$ has no loops. An 
absolute upper bound for $b_i$ when $1 \le i \le r$ is  
$\binom{n-i-1}{r-i}$ and this is achieved by $U_{n,r}.$
Theorem \ref{CM h-vectors} gives a  relative upper bound of 
$b_{r-i} \le b_{r-i+1}^{<i-1>} .$
Absolute lower bounds for $b_i$ were determined by Brylawski.

\begin{thm} \cite{Bry3}
  If $M$ is as above, then $b_i \ge n-r$ for all $i, 2 \le i \le r-1.$  
\end{thm}

In order to find relative lower bounds for  $b_1$ we introduce the 
following definition.

\begin{defn}
  Let $S$ be a series class of a connected matroid $M.$
  Then $S$ is a {\bf regular series} class of $M$ if $M-S$ is connected.  
\end{defn}

\begin{prop}
  If $M$ is connected and contains more than one series class, then $M$ contains at least three regular series classes.
\end{prop}
\begin{proof}
  Induction on  $m$, the number of series classes in $M.$  A matroid 
  with exactly two series class is not connected.  If $m=3,$ then $M$ is 
  the cycle matroid of a theta graph with exactly three paths. In this case all three of the series classes are regular.   

  For the induction step, let $S$ be a series class which is not regular.  Let $\tilde{M}$ be the matroid obtained by contracting
  all but one of the elements of $S.$ Let $e$ be the remaining element of $S.$
  Since $\tilde{M}$ is connected, but $\tilde{M}-e$ is not connected, $\tilde{M}$  
  is the  series connection of two connected matroids $A$ and $B$ at $e$ \cite[Theorem 7.1.16]{O}.  Both $A$ and $B$ must contain more than one series class, otherwise they would be contained in $S.$  Therefore, the induction 
  hypothesis applies to $A$ and $B.$  Even if $\{e\}$ is contained in 
  a regular series class in $A$ and $B,$ both $A$ and $B$ contain two other regular series classes.  All four of these series classes are regular in $M.$
\end{proof}

\begin{thm} \label{char by r}
  If $M$ is connected and $1 \le i \le r$, then
\begin{equation} \label{eq1}
b_i \le \binom{r-2}{i-1} b_1 + \binom{r-2}{i-2}.
\end{equation}
\end{thm}

\begin{proof}
  The proof is by   
  induction on $n,$ the initial case being the three--point line.    Let $S$ be a series class of $M.$  If $S$ is the only series class of $M,$ then $M$ is a circuit and (\ref{eq1}) holds. Otherwise,  by the previous proposition, we may choose  $S$ to be a regular series class. In particular, $M-S$ is connected.  We break the induction step into three  cases.

\begin{enumerate}
  \item
     $M-S$ and $M/S$ are connected:
  Let $s = |S|.$ If $s > i,$
  then $b_i(M) = 
  b_i(\hat{M})$ and $b_1(M) = b_1(\hat{M}),$
  where $\hat{M}$ is $M$ with $S$ contracted down to a series class 
  of cardinality $i.$  So, we will assume that $s \le i.$   Let $\tM$ be $M$ 
  with $S$ contracted down to a single element $e.$  Since $M$ is connected, $e$ is neither a loop nor a coloop of $M.$  Applying Tutte recursion   to $M$ and then again to $\tM$ we see that

$$b_i(M) = b_i(\tM/e) + \sum^{s-1}_{j=0} b_{i-j}(\tM-e).$$
 
 Now, since $\tM/e = M/S$ is a rank $r-s$ connected matroid and $\tM-e = M - S$ is a rank $r-s+1$ connected matroid, the induction hypothesis implies that the above expression is bounded above by
$$\binom{r-s-2}{i-1} b_1(\tM/e) + \binom{r-s-1}{i-2} +  
\sum^{s-1}_{j=0} \binom{r-s-1}{i-j-1} 
b_1(\tM-e) + \sum^{s-1}_{j=0} \binom{r-s-1}{i-j-2} $$

$$\le \binom{r-2}{i-1} b_1(\tM/e) +  \binom{r-2}{i-1} b_1 (\tM-e) + 
\binom{r-2}{i-2} +$$
$$\left\{ \binom{r-s-2}{i-1} - 
\binom{r-2}{i-1} \right\} b_1(\tM/e) + \binom{r-s-2}{i-2}.$$

\noindent Since $\tM/e$ is connected, $b_1(\tM/e) \ge 1.$ 
Thus, the last row is non--positive and (\ref{eq1}) is satisfied.  To see the last inequality, note that 
 $$\sum^{s-1}_{j=0} \binom{r-s-1}{i-j-1} \le \sum^{s-1}_{j=0} \binom{r-s-1}{i-j-1} \binom{s-1}{j}= \binom{r-2}{i-1},$$

 \noindent and similarly,
 
 $$\sum^{s-1}_{j=0} \binom{r-s-1}{i-j-2} \le \sum^{s-1}_{j=0} \binom{r-s-1}{i-j-2} \binom{s-1}{j}= \binom{r-2}{i-2}.$$

\item
  $S = \{e\}, M-e$ is connected, but $M/e$ is not connected:
    Then, $M$ is the parallel connection of two connected 
  matroids $A$ and $B$ with $r(A) + r(B) -1 = r$ \cite[Theorem 7.1.16]{O}.  By Theorem \ref{tuttepoly}
  and the 
  induction hypothesis,

$$b_i(M) = \sum_{j+k-1=i} b_j(A) b_k(B)$$
$$ \begin{array}{ll}
\le & \displaystyle\sum_{j+k-1=i} \left\{ \binom{r(A) 
-2}{j-1} 
b_1(A) + \binom{r(A)-2}{j-2} \right\} \left\{ \binom{r(B) -2}{k-1} 
b_1(B) + \binom{r(B)-2}{k-2} \right\} \\ 
\ & \ \\
 = &  \displaystyle\sum_{j+k-1=i} \left\{ \binom{r(A)-2}{j-1} \binom{r(B)-2}{k-1} b_1(A) 
b_1(B) + \binom{r(A)-2}{j-1} \binom{r(B)-2}{k-2} b_1(A) \right\} + \\
\ & \\
\ & \displaystyle\sum_{j+k-1=i} \left\{  \binom{r(A)-2}{j-2} \binom{r(B)-2}{k-1} b_1(B) + 
\binom{r(A)-2}{j-2} 
\binom{r(B)-2}{k-2} \right\} \\
\ & \\
= & \displaystyle\binom{r-3}{i-1} b_1(A) b_1(B) + \binom{r-3}{i-2} 
b_1(A) + \binom{r-3}{i-2} b_1(B) + \binom{r-3}{i-3}.
\end{array}$$ \\

\noindent Therefore,

$$\begin{array}{ll}
\ & \displaystyle\binom{r-2}{i-1} b_1(M) + \binom{r-2}{i-2} - b_i(M) \\
\ge &  \displaystyle\left\{ \binom{r-2}{i-1} - \binom{r-3}{i-1} \right\} b_1(A) b_1(B) + 
\left\{ \binom{r-2}{i-2} - \binom{r-3}{i-3} \right\} - \\
\ & \displaystyle\binom{r-3}{i-2} \left\{ b_1(A) + b_1(B) \right\} \\
= & \displaystyle\binom{r-3}{i-2} \left( b_1(A) b_1(B) + 1 - b_1(A) - b_1(B) \right) 
\ge 0.
\end{array}$$

\item
Finally, suppose that $S$ is a non--trivial series, $M-S$ is connected, but $M/S$ is not connected.  Let $s, \tM$ and $e$ 
be as above. Since $\tM/e$ is not connected, $b_1(\tM)=b_1(\tM-e).$  
Therefore,

$$b_i(M) = b_i(\tM) + \sum^{s-1}_{j=1} b_{i-j}(\tM-e) $$
$$\le b_1(M) \left\{ \sum^{s-1}_{j=0} \binom{r-s-1}{i-j-1} \right\} + 
\sum^{s-1}_{j=0} \binom{r-s-1}{i-j-2} $$
$$ \le \binom{r-2}{i-1} b_1(M) + \binom{r-2}{i-2}.$$

\end{enumerate}
\end{proof}

\begin{cor} \label{disconnected}
  Let $M$ be a rank $r$ matroid with $k$ components, $r-k \ge 2.$  
  Let $2 \le i \le r-k.$  Then,
\begin{equation} \label{eqdisc}
b_{i+k-1}(M) \le \binom{r-k-1}{i-1} b_k(M) + \binom{r-k-1}{i-2}.
\end{equation}
\end{cor}

\begin{proof}
 Since $k=1$ is the previous theorem we assume that $M$ is not connected.
 Let $M=M_1 \oplus \dots \oplus 
M_k$ be a direct sum decomposition of $M$ into connected matroids.  
Define $\tM_1=M_1.$  Given $\tM_i$ let $\tM_{i+1}$ be any parallel 
connection of $\tM_i$ and $M_{i+1}.$  Then $\tM_k$ is a connected matroid 
of rank $r-k+1.$  Furthermore, by Theorem 
\ref{tuttepoly} $b_{i+k-1}(M) = b_i(\tM_k).$
Since (\ref{eq1}) holds for the connected $\tM,$ (\ref{eqdisc}) 
holds for $M.$   
\end{proof}

When does equality occur in the above theorem?  The proof 
shows that if equality occurs, then it must also occur in the minors of 
$M$ used in the induction.  Combining this with an induction argument shows 
that if $b_i(M) = \binom{r-2}{i-1} b_1(M) + \binom{r-2}{i-2},$ then $b_j(M) = 
\binom{r-2}{j-1} b_1(M) + \binom{r-2}{j-2}$ for all $1 \le j \le i.$  
Brylawski  proved (\ref{eq1}) for $i=r-1.$  
He also showed that given $b_1$ and $r,$ then equality occurs if $M$ is 
the parallel connection of a $(b_1+2)$--point 
line and $r-1$ three--point lines.  Hence, (\ref{eq1}) is optimal, although a complete description of 
the matroids 
which satisfy equality in this corollary remains unknown \cite{Bry3}.  

The coefficient $b_1(M)$ is also known as $\beta(M),$ the {\it beta} invariant of 
$M.$  Brylawski identified matroids with beta invariant 1 as 
series--parallel matroids \cite{Bry5} while Oxley classified matroids with 
$2 \le \beta(M) \le 4$ \cite{Ox2}.

\begin{thm}
  Assume $r \ge 2$ and let $\beta=b_1(M).$ Then, for all $i, 0 \le i \le r,$

$$w_i \le \sum^i_{j=0} \binom{r-j}{r-i} \left\{ \binom{r-2}{r-i-1} \beta + 
\binom{r-2}{r-i-2} \right\}.$$
\end{thm}

\begin{proof}  
This follows immediately from (\ref{h by f}) and Theorem \ref{char by r}.
\end{proof}

It is also possible to estimate  $b_i$ in terms of $n-r.$  
For positive integers $i$ and $x$  define 
$$\phi_i(x)=\binom{x-2}{i-1} \binom{x-1}{0} + \binom{x-2}{i-2} 
\binom{x}{1}  + \dots + \binom{x-2}{0} 
\binom{x+i-2}{i-1}.$$

\begin{thm}
  Suppose $M$ is connected.   Then,
\begin{equation} \label{eq2}
  b_i(M) \le \phi_i(n-r) b_1(M) + \phi_{i-1}(n-r).
\end{equation}
\end{thm}
  
\begin{proof}
  We can assume that every series class of $M$ has exactly $i$ elements.  Indeed, by [a] and [c] of Theorem \ref{tuttepoly}, any series class with more than $i$ elements can be contracted down to cardinality $i$ without changing either side of (\ref{eq2}), while expanding any class with fewer than $i$ elements may increase the left--hand side of (\ref{eq2}) but will not alter the right--hand side.  Let $\tM$ be the matroid obtained from $M$ by contracting all of the series classes down to one element.  The dual of the formula on the top of page 185 of \cite{Bry2} is

\begin{equation} \label{eq3}
T(M;x,0) = (x^{i-1} + \dots x +1)^{n-r} T(\tM; x^i, \frac{x^{i-1} + \dots x}{x^{i-1} + \dots x + 1})
\end{equation}

 Using (\ref{eq3}), we see that,

\begin{equation} \label{eq4}
  b_i(M) = \displaystyle\sum^i_{j=1} \binom{n-r+i-j-1}{i-j} b^\star_j(\tM).
\end{equation}

\noindent  Since $b^\star_1(\tM) = b_1(M),$ (\ref{eq2}) follows from (\ref{eq4}) by applying (\ref{eq1}) to $\tM^\star.$
   
\end{proof}
   
Inequality (\ref{eq2}) is as optimal as can be expected in the sense that given $n-r,i$
and $b_1$ there are matroids which satify equality.  Take any matroid which 
satisfies equality in (\ref{eq1}) and expand every series class to 
cardinality $i.$ Then, equality in (\ref{eq2}) holds.
Of course, since $b_r=1$ and $\phi_i$ is increasing in 
$i$, no matroid can satisfy equality in (\ref{eq2}) for all $i.$

\section{Independence complexes} \label{ind complexes}

Suppose the  smallest cocircuit of $M$ has 
cardinality $k.$   As 
pointed out in Section \ref{k-CM}, $\Delta(M)$ is a $k$--CM complex.  So, we can 
apply those methods to $\Delta(M).$ In addition to the previously 
mentioned absolute upper bound $h_i(M) \le \binom{n-r+i-1}{i}$ and relative 
upper bound $h_{i+1} \le h^{<i>}_i,$ the $h$--vectors of independence 
complexes of matroids satisfy an analogue of the $g$--theorem for simplicial 
polytopes.

\begin{thm} \cite{Sw3}
Assume that $M$ has no  coloops.  Let 
$g_i(M)=h_i(M)-h_{i-1}(M).$  Then for all $i, 1 \le i \le (r+1)/2,$
$$g_{i+1}(M) \le g_i^{<i>}(M).$$
\end{thm}

\noindent The above theorem was proved independently by Hausel and Sturmfels for matroids representable over the rationals using toric 
hyperk\"{a}ler varieties \cite{HS}.

Relative lower bounds, also reminiscent of the $g$--threorem for 
simplicial polytopes, were originally established by Chari using a PS--ear 
decomposition of $\Delta(M).$  See \cite{Ch} for the definition of PS--ear 
decompositions and a proof of the following theorem.

\begin{thm}
  Suppose $M$ has no coloops.  Then for all $i, 0 \le 
  i \le r/2.$
$$ h_{i-1} \le h_i,$$
$$ h_i \le h_{r-i}.$$
\end{thm}

\begin{prob}
  Do $2$--CM complexes satisfy the inequalities in the previous two 
  theorems?
\end{prob}

An affirmative answer to this question would, with the addition of the 
Dehn--Sommerville equations, give a complete description of all possible 
$h$--vectors of simplicial homology spheres \cite[Conjecture II.6.2]{St}.

In \cite{BC2} Brown and Colbourn conjectured that for co-graphic $M,$ the complex zeros of $T(M;x,1)$ were contained in the closed unit disk.  While this has since proven to be false \cite{RS}, attempts to prove it led to a couple of relative lower bounds for $h$--vectors of independence complexes of any matroid.

\begin{thm}  Suppose $M$ has no coloops.
 \begin{enumerate}
   \item \cite{BC2} For all $i \le r,$
       $$ h_i \ge \sum^i_{j=1} (-1)^{j-1} h_{i-j}$$
       
   \item \cite{Wa}
       Let $ I_j$ be the number of independent subsets of $M$ of cardinality $j.$ Then for all $0 \le k \le r,$ 
$$ \sum^r_{j=k} \binom{j}{k} (-2)^{r-j} I_j \ge 0.$$
 \end{enumerate}
\end{thm}

Stanley used the notion of a level ring to establish the relative lower 
bound $ 
h_{j-i}(M) \le h_i(M) h_j(M)$ whenever 
$0 \le i,j \le r.$  In particular, setting $j=r,$ we find that $h_{r-i}(M)\le 
\binom{n-r+i-1}{i} h_r(M).$  By applying (\ref{del by h}) we can 
obtain similar relative lower bounds for $h_{i-j}(M)$ in terms of $h_i(M)$ 
and we can also determine when equality occurs.

\begin{prop}
  Assume that $M$ has no coloops.    Then 
  for all $i, 1 \le j < i \le r,$
\begin{equation} \label{max h}
  h_{i-j}(M) \le \frac{\binom{n-i+j-1}{r-i+j}}{\binom{n-i-1}{r-i}} h_i (M).
\end{equation}
 
Furthermore,  equality occurs if and only if every series class of $M$ has 
cardinality greater than $r-i+j.$
\end{prop}

\begin{proof}  
Since $M$ has no coloops, $\Delta(M)$ is a $2$--CM 
complex.  Therefore,  (\ref{del by h}) implies $(r-i+1) h_{i-1}(M) \le (n-i) 
h_i(M).$  In order for equality to occur, $h_i(M-e)$ must be zero for 
every $e$ in $E.$ By Theorem \ref{tuttepoly} [a], this is equivalent to every 
series class of $M$ having cardinality greater than  $r-i+1.$  The 
proposition 
follows by induction on $j.$
\end{proof}
 
 In \cite{BC} Brown and Colbourn proved the relative lower bound $h_{r-1} (M) 
\le r h_r(M)$ which only involves the rank of $M.$  This can be improved 
using Theorem \ref{char by r}.  

\begin{thm}
  Let $M$ be a rank $r$ matroid without coloops.  Then,
\begin{equation} \label{ind by r}
h_{r-i} \le \binom{r-1}{i} h_r + \binom{r-1}{i-1}.
\end{equation}
\end{thm}

\begin{proof}
  By \cite{Bry4}, $h_i(M)$ equals $b_{r-i+1}$ of the free coextension of 
  $M.$  Since the latter matroid has rank $r+1$ and is connected, 
  (\ref{ind by r})  is an immediate consequence of Theorem \ref{char by r}.
\end{proof}

As in the case of Theorem \ref{char by r}, if $h_{r-i}(M) = 
\binom{r-1}{i} h_r(M)+ \binom{r-1}{i-1},$ then $h_{r-j}(M) \le \binom{r-1}{j} 
h_r(M) + \binom{r-1}{j-1}$ for all 
$0 \le j \le i.$ 
A routine 
deletion--contraction induction shows that for a given $r$ and $h_r,$
$$M=U_{1,h_r+1} \oplus \underbrace{U_{1,2} 
\oplus \dots \oplus U_{1,2}}_{r-1}$$
satisfies equality in (\ref{ind by r}).  

\begin{cor}
  Let $M$ be a rank $r$ matroid without coloops. Let $I_j$ be the number 
  of independent subsets of $M$ of cardinality $j.$ Then,
$$I_j \le \displaystyle\sum^j_{i=0} \binom{r-i}{r-j} \left\{ 
\binom{r-1}{i} h_r + \binom{r-1}{i-1} \right\}.$$
\end{cor}

\begin{proof}  Apply the above theorem to (\ref{f by h}). \end{proof}

In section \ref{k-CM} we posed the problem of finding absolute 
lower bounds for a $k$--CM complex given $n$ and $r.$ Here we examine this problem for independence complexes. Consider the special case of a rank two matroid $M$ without loops.   The simplification of $M$  is isomorphic to $U_{2,m}$  where $m$ is the number of parallel classes of $M.$    Therefore, $M$ is specified up to isomorphism by a partition $n = p_1 + \dots + p_m,$ where the $p_i$'s are the sizes of the parallel classes of $M.$   Since $h_0 = 1$ and $h_1 = n-r,$ minimizing the $h$--vector of $M$ is equivalent to minimizing the number of bases of $M.$  As noted earlier, $M$ is $k$--CM if and only if every hyperplane of $M$ has cardinality at most $n-k.$  Equivalently, each $p_i \le n-k.$  The number of bases of $M$ is 

$$\binom{n}{2} - \displaystyle\sum^m_{i=1} \binom{p_i}{2}.$$

\noindent This is minimized by setting $m = \lceil n/(n-k) \rceil, p_i = n-k$ for $i \le m-1,$ and $p_m = n - (m-1)(n-k).$  Note that this implies that when $n \ge 2k, h_2(M)$ is bounded below by $h_2(U_{1,n-k} \oplus U_{1,k}).$

 An independence 
complexes is $2$--CM if and only if it has no coloops.  In \cite{Bj} Bj\"{o}rner showed that for any matroid without coloops $h_i \ge n-r$ for $0 < i < r.$  While it is not specifically stated, the proof implies that $h_r \ge n-2r+1.$  In general, given $n$ and $r$ there may be no single coloop-free matroid that achieves all of these bounds.  For example, if $n=8$ and $ r=4,$  then the only  matroid without coloops such that  $h_4(M) = 1$ is $M =U_{1,2} \oplus U_{1,2} \oplus U_{1,2} \oplus U_{1,2}.$  However, $h_2(M) = 6 > n-r.$  If we restrict our attention to $i < r, $ then $U_{1,n-r} \oplus U_{r-1,r}$ does satisfy $h_i = n-r$ for $0 < i < r.$

\begin{defn}
  $M(r,n,k) =U_{1,n-r-k+2} \oplus U_{r-1,r+k-2}$
\end{defn} 

\noindent Direct computation shows that $h_i(M(r,n,k)) = \binom{k+i-2}{i} + (n-r-k+1) \binom{k+i-3}{i-1}.$ In addition, $\Delta(M(r,n,k))$ is $k$--CM as long as $n \ge r+2k-2.$

\begin{thm} \label{matroid k-cm}
  Fix $r \ge 2$ and $k \ge 3$ .  There exists $N(k,r)$ such that if $M$ is a matroid without loops whose smallest cocircuit has cardinality at least $k$ and $n \ge N(k,r),$ then for all $i, 0 \le i \le r,$

\begin{equation} \label{long term}
h_i(M) \ge h_i(M(r,n,k)).
\end{equation}

\end{thm}

\begin{proof}
  First we show that if $n > k(r+1),$ then there exists $e \in M$ such that $\Delta(M-e)$ is still $k$--CM.  Let $\mathbf{H}$ be the set of hyperplanes of $M$ of cardinality $n-k.$  If $\mathbf{H}$ is empty, then any $e$ will do since no hyperplane of $M-e$ will have size greater than $n-k-1.$  Otherwise, let $B$ be the intersection of all of the hyperplanes in $\mathbf{H}.$  Since $B$ is a flat of $M$ there exists $H_1, \dots, H_{r+1},$ not necessarily distinct, in $\mathbf{H}$ such that $H_1 \cap \dots \cap H_{r+1} = B.$  Therefore, $|B| \ge n-k(r+1)$ and $B$ is not empty.  But, for any $e \in B, \Delta(M-e)$ is $k$--CM.

As noted above, when $r=2, N(2,k) = 2k$ works.  So,  assume that $r \ge 3.$  Let $M^\prime$ be a contraction of $M$ and let $n^\prime=|E(M^\prime)|.$  By Proposition \ref{k-CM abs},  $h_i(M^\prime) \ge h_i(U_{r-1,r+k-2}).$  In fact, if $n > r+k-2,$ then $h_i(M^\prime)$ is strictly greater than $h_i(U_{r-1,r+k-2})$ for $1 \le i \le r.$ Indeed,  this is proved by Tutte recursion as in Proposition \ref{k-CM abs}.  The base case compares the $h$--vectors of $U_{2,4}$ and any five element rank two matroid whose smallest cocircuit has at least three elements.  The $h$--vector of $U_{2,4}$ is $(1,2,3).$ From the discussion of rank two matroids, the $h$--vectors of the latter group of matroids is bounded below by $(1,3,4),$  the $h$--vector of the matroid whose simplification is $U_{2,3}$ and whose parallel classes have cardinality $2,2$ and $1.$  Note that this claim is not true when $k=2.$  In particular, $U_{1,2} \oplus U_{r-1,r}$ is a rank $r$  matroid without coloops and $r+2$ elements whose $h_r$ is not strictly less than $h_r$ of $U_{r,r+1}.$  

To finish the proof, we find $N(r,k,i)$ such that the theorem holds for just $h_i$ and then let $N(r,k)$ be the maximum of the all of the $N(r,k,i).$ Since $h_0(M)=1$ and $h_1(M) = n-r, r+k-1$ works for $N(r,k,0)$ and $N(r,k,1).$ So fix $i \ge 2.$ Let $N$ be the minimum of $h_i(\bar{M})$ for all loopless matroids $\bar{M}$ such that $|E(\bar{M})| = k(r+1) +1, r(\bar{M}) = r$ and the smallest cocircuit of $\bar{M}$ has at least $k$ elements. Let $N(r,k,i) = k(r+1) + 1 + h_i(M(r,k(r+1)+1,k)) - N.$ 
\\ 

{\bf Claim:} If $n \ge N(k,r,i),$ then $h_i(M) \ge h_i(M(r,n,k)).$ \\

{\it Proof of claim:} Choose $e_1 \in M$ such that the smallest cocircuit of $M-e_1$ has cardinality greater than or equal to $k.$  Given $e_j$ choose $e_{j+1}$ so that the smallest cocircuit of $M - \{e_1,\dots,e_j,e_{j+1}\}$ has size at least $k.$  This can be done up to $j = n-k(r+1)-1.$  Deleting and contracting on each deletion,

$$h_i(M) = h_i(\tM) + \displaystyle\sum_j h_{i-1}(M - \{e_1,\dots,e_{j-1}\}/e_j),$$

\noindent where $\tM$ is $M - \{e_1,\dots,e_{n-k(r+1)-1}\}.$  By construction, $|E(\tM)| = k(r+1)+1, r(\tM) = r$ and the smallest cocircuit of $\tM$ has at least $k$ elements.  In addition, the rank of each contraction is $r-1$ and its independence complex is $k$--CM.  There are two possibilities.

\begin{itemize} 

 \item
    Every contraction has more than $r+k-2$ non-loop elements. In this case $h_i(M) \ge h_i(\tM) + (n-k(r+1)-1)[ h_{i-1}(U_{r-1,r+k-2})+1].$ Compare this to computing $h_i(M(r,n,k))$ by deleting and contracting down to $U_{1,rk-k+3} \oplus U_{r-1,r+k-2}.$ The definition of $N(r,k,i)$ insures that $h_i(M)$  is bounded below by $h_i(M(r,n,k)).$
    
 \item
    At least one contraction, say $M - \{e_1,\dots,e_{j-1}\}/e_j)$ has exactly $r+k-2$ elements.  Since this contraction is a rank $r-1$ matroid whose smallest cocircuit has at  least  $k$ elements it must be equal to $U_{r-1,r+k-2}.$  Therefore, $M - \{e_1,\dots,e_{j-1}\}$ has one non--trivial parallel class which contains $e_j$ and the simplification of $M - \{e_1,\dots,e_{j-1}\}$ is a one--element coextension of $U_{r-1,r+k-2}.$   The one--element coextension of $U_{r-1,r+k-2}$ which minimizes $h_i(M - \{e_1,\dots,e_{j-1}\})$ is the one obtained by adding a coloop to $U_{r-1,r+k-2}.$  Hence, $h_i(M - \{e_1,\dots,e_{j-1}\})$ is bounded below by $h_i(M(r,n-j,k)).$ However,  this implies that $h_i(M) \ge h_i(M(r,n-r,k) + j \ h_{i-1}( U_{r-1,r+k-2}) = h_i(M(r,n,k)).$
    
\end{itemize}

\end{proof}

Some lower bound on $n$ is necessary in order for (\ref{long term}) to hold. For instance, let $M= U_{1,3} \oplus U_{1,3} \oplus U_{1,3}.$  Then $r=3, k=3$ and $n=9.$  The $h$--vector of $M$ is $(1,6,12,8),$ while the $h$--vector of $M(3,9,3) = U_{1,5} \oplus U_{2,4}$ is $(1,6,11,12).$

\noindent As usual, absolute lower bounds yield relative lower bounds via (\ref{del by h}).
\begin{cor}
 Fix $r \ge 2$ and $k \ge 3$ .  There exists $N(k,r)$ such that if $M$ is a matroid without loops whose smallest cocircuit has cardinality $k$ and $n \ge N(k,r),$ then for all $i, 0 \le i \le r,$

$$(r-i+1) h_{i-1}(M) + n \ h_i(M(r,n-1,k-1)) \le (n-i) 
h_i(M).$$
\end{cor}

\noindent {\it Acknowledgment:} An anonymous referee's comments and suggestions dramatically improved the exposition in several places.

 \end{document}